\documentclass[11pt]{article}
\usepackage[dvips]{graphicx}
\usepackage{amsfonts,epsf,amsmath,tikz}

\newtheorem{thm}{Theorem}
\newtheorem{lem}[thm]{Lemma}

\newtheorem{ques}[thm]{Question}

\newcommand{\proof}{\noindent{\bf Proof.} }
\newcommand{\cp}{\,\square\,}
\newcommand{\qed}{\begin{flushright} $\Box$ \end{flushright}}
\begin{document}

\title{On the Cartesian product of non well-covered graphs}

\author{Bert L. Hartnell\thanks{This author is partially supported
by NSERC.}
\\
        Saint Mary's University\\
        Halifax, Nova Scotia, Canada\\
        bert.hartnell@smu.ca
\and
Douglas F. Rall\thanks{The second author is Herman N. Hipp Professor of Mathematics
at Furman University.
This work was partially supported by a grant from the Simons Foundation
(\#209654 to Douglas Rall).
Part of the research done during a sabbatical visit at Saint Mary's University.}
\\
Furman University\\
Greenville, SC, USA\\
doug.rall@furman.edu
}

\date{\today}

\maketitle

\begin{abstract}
A graph is well-covered if every maximal independent set has the same
cardinality, namely the vertex independence number.  We answer a question
of Topp and Volkmann~\cite{tv1992} and prove that if the Cartesian product
of two graphs is well-covered, then at least one of them must be well-covered.
\end{abstract}

{\small {\bf Keywords:} maximal independent set, well-covered, Cartesian product} \\
\indent {\small {\bf AMS subject classification: 05C69, 05C76}}

\section{Introduction}

 A {\it well-covered graph} $G$ (Plummer~\cite{mdp1970}) is one in which every maximal independent
set of vertices has the same cardinality. That is, every maximal independent set (equivalently,
every independent dominating set) is a maximum independent set.
This class of graphs has been investigated by many researchers from several different
points of view.  Among these are attempts to characterize those well-covered graphs with a
girth  or a maximum degree restriction.  For more details on these approaches as well as
others see the surveys by Plummer~\cite{mdp1993} and by Hartnell~\cite{h1999}.

Topp and Volkmann~\cite{tv1992} investigated how the standard graph
products interact with the class of well-covered graphs.  They asked the following
question that was restated by Fradkin~\cite{f2009}.

\begin{ques} {\rm (\cite{tv1992})} \label{ques:1}
Do there exist non well-covered graphs whose Cartesian product is well-covered?
\end{ques}

The principal result of this paper is the following theorem that answers Question~\ref{ques:1} in the
negative.

\begin{thm} \label{thm:main}
If $G$ and $H$ are graphs whose Cartesian product is well-covered, then at least one
of $G$ or $H$ is well-covered.
\end{thm}

In Section~\ref{sec:defs} we define the terms used most often in this paper; standard
graph theory terminology is used throughout. We then establish Theorem~\ref{thm:main}
in Section~\ref{sec:mainresults}.

\section{Definitions} \label{sec:defs}

If $G_1=(V_1,E_1)$ and $G_2=(V_2,E_2)$ are any two graphs, the {\it Cartesian
product} of $G_1$ and $G_2$ is the graph denoted $G_1\cp G_2$ whose vertex
set is the Cartesian product of their vertex sets $V_1 \times V_2$.  Two vertices
$(x_1,x_2)$ and $(y_1,y_2)$ are adjacent in $G_1\cp G_2$ if either
$x_1=y_1$ and $x_2y_2\in E_2$, or $x_1y_1 \in E_1$ and $x_2=y_2$. Note that if
$I_1$ is independent in $G_1$ and $I_2$ is independent in $G_2$, then
the set $I_1\times I_2$ is independent in $G_1\cp G_2$.

For an arbitrary graph $G$ we follow Fradkin~\cite{f2009} and define
a {\it greedy independent decomposition}  of $G$ to be a partition
$A_1,A_2, \ldots, A_t$ of $V(G)$ such that $A_1$ is a maximal independent set in $G$, and
 for each $2 \le i \le t$, the set $A_i$ is a maximal independent set in the
graph $G-(A_1\cup \cdots\cup A_{i-1})$.
One way to construct maximal independent sets in the Cartesian product
$G\cp H$ is to select any greedy independent decomposition $A_1,A_2, \ldots, A_t$
of $G$ and an arbitrary greedy independent decomposition $B_1,B_2, \ldots, B_s$ of
$H$ and combine them into what is called a ``diagonal'' set of the product
as $M=\cup_i(A_i \times B_i)$.  If $s \neq t$, then there are as many sets in this
union as the smaller of $s$ and $t$.

The {\it vertex independence number} of a graph $G$ is the cardinality of a largest
independent set in $G$.  We denote the vertex independence number of $G$ by $\alpha(G)$
and refer to an independent set of this order as an $\alpha(G)$-set.
If a graph $G$ has an independent set $M$ such that $G-N[M]=\{x\}$
for some vertex $x$, then $x$ is said to be an {\em isolatable vertex} of
$G$.  The existence of such a vertex is central to our work.

\begin{lem} \label{lem:cliqueremainder}
Let $G$ be a graph in which no vertex is isolatable.  If $I$ is any maximum
independent set in $G$ and $x$ is any vertex of $I$, $G-N[I-\{x\}]$ is a clique
of order at least two.
\end{lem}
\proof
Suppose $I$ is an $\alpha(G)$-set and that $I$ has a vertex $v$ such that the graph
$G-N[I-\{v\}]$ has an independent set $A$ of size at least two.  Then
$I'=(I -\{v\})\cup A$ is independent in $G$ and has order larger than $|I|=\alpha(G)$,
a contradiction.  \qed

\section{Main Results} \label{sec:mainresults}

We first reduce the study of when a Cartesian product is well-covered
by considering the existence of isolatable vertices in the two factors.

\begin{thm} \label{thm:firstone}
Suppose that $H$ is not well-covered and $G$ has an isolatable vertex.
Then $G \cp H$ is not well-covered.
\end{thm}
\proof
Let $A$ and $B$ be maximal independent subsets of $H$ with $|A|>|B|$,
and suppose that $x$ is an isolatable vertex in $G$.  Let $I$ be an
independent set in $G$ such  that $x$ is an isolated vertex in the graph
$G-N[I]$.  Extend the independent
set $I \times A$ to a maximal independent set $J$ of $(G-N[x])\cp H$.
Let $m=|J|$.  Note that $J$ dominates $N_G(x) \times A$ (and perhaps
other vertices of $N_G(x) \times V(H)$), but $J$ does not contain any
vertices from $N_G[x] \times V(H)$.

Let $J_1=J \cup (\{x\} \times A)$ and $J_2=J \cup (\{x\} \times B)$.
By the choice of $A$ and $B$ it is clear that $|J_1| > |J_2|$.  Let
$X_A$ denote the set of vertices in $N_G(x) \times V(H)$ that are not
dominated by $J_1$.  Similarly, let $X_B$ denote the set of vertices
in $N_G(x) \times V(H)$ that are not dominated by $J_2$.  The set
$X_B$ is a subset of $X_A$.

Choose a maximal independent set $L$ of the subgraph of $G \cp H$
induced by $X_B$.  Then $J_2 \cup L$ is a maximal independent set
in $G \cp H$.  Extend $L$ to a maximal independent set $M$ of the
subgraph of $G \cp H$ induced by $X_A$.  Now, $J_1 \cup M$ is a
maximal independent set of $G \cp H$, and
\[|J_1 \cup M| = |J_1|+|M|>|J_2|+|M| \ge |J_2|+|L|=|J_2 \cup L|\,.\]
Therefore, $G\cp H$ has maximal independent sets of distinct
cardinalities, and thus $G \cp H$ is not well-covered. \qed

It now follows that if both of $G$ and $H$ are not well-covered but $G\cp H$ is
well-covered, then neither $G$ nor $H$ has an isolatable vertex.

\begin{lem} \label{lem:samesizedisjoint}
Let $G$  and $H$ be graphs such that neither has an isolatable vertex.
If $G\cp H$ is well-covered, then  both
$G$ and $H$ have the property that if $M$ is any maximal independent set
of the graph, that graph must have a maximal independent set $N$ that is
disjoint from $M$.  Furthermore, at least one of $G$ or $H$ has the property
that any two disjoint maximal independent sets have the same cardinality.
\end{lem}
\proof
As observed above, neither $G$ nor $H$ has an isolatable vertex.
Let $I_1,I_2, \ldots, I_t$ be  any greedy independent decomposition of $G$ and
$J_1,J_2, \ldots, J_s$ be any greedy independent decomposition of $H$.  Let $p=\min\{s,t\}$.
The (so-called ``diagonal'') set
$M= (I_1\times J_1) \cup (I_2\times J_2) \cup \cdots \cup (I_p\times J_p)$
is maximal independent in $G\cp H$.  Since $G \cp H$ is well-covered, this
implies that $\alpha(G\cp H)=|M|=\sum_{k=1}^p |I_k|\!\cdot\!|J_k|$.

Since $J_1$ is an independent set in $H$ and $J_2$ is a maximal independent
set in $H-J_1$, if there exists a vertex $u\in J_1$ that is not dominated by
$J_2$, then $u$ is isolated in $H-N[J_2]$.  This contradicts the fact that
$H$ does not have an isolatable vertex.  Therefore, $J_2$ is actually a maximal
independent set in $H$ as well as in $H-J_1$.  By an identical argument it follows that
$I_2$ is a maximal independent set in $G$.  Suppose that $a=|J_1|$ and $b=|J_2|$
and that $a\neq b$.  Let $c=|I_1|$ and $d=|I_2|$.  Since $I_1$ and $I_2$ are
disjoint maximal independent sets in $G$, the list
$I_2,I_1,I_3,\ldots,I_t$ is also a greedy independent decomposition of $G$.  This implies
\[ca+db+\sum_{k=3}^p|I_k|\!\cdot\!|J_k| =\alpha(G \cp H)=da+cb + \sum_{k=3}^p|I_k|\!\cdot\!|J_k|\,,\]
since $G\cp H$ is well-covered, and thus $ca+db=da+cb$.  Since $a\neq b$ we get $c=d$; that is, $|I_1|=|I_2|$.
Since $I_1,I_2, \ldots, I_t$ is an arbitrary greedy independent decomposition of
$G$, the lemma follows. \qed

We now proceed to prove our main result.
\setcounter{thm}{1}
\begin{thm} \label{thm:main}
If $G$ and $H$ are graphs such that $G\cp H$ is well-covered, then at least one
of $G$ or $H$ is well-covered.
\end{thm}
\proof
Suppose by way of contradiction that the statement is not true.  Let $G$
and $H$ be a pair of graphs neither of which is well-covered but such that
$G \cp H$ is well-covered.  As above we may assume that no vertex of
either $G$ or $H$ can be isolated in its own graph.  From Lemma~\ref{lem:samesizedisjoint}
we may assume without loss of generality that $G$ has the property that any
two maximal independent sets of different cardinalities must intersect nontrivially.

Since $G$ is not well-covered, there exists a maximal independent set whose cardinality
is less than $\alpha(G)$.  From the collection of all maximal independent sets in $G$
choose a pair, say $I$ and $J$,  such that $|J|<|I|=\alpha(G)$ and
$|I \cap J|$ is as small as possible.
Since $|I| \neq |J|$ there exists $v \in I\cap J$.  Let $F=G-N[I-\{v\}]$.
By Lemma~\ref{lem:cliqueremainder} this subgraph $F$ is a clique of order at least two.
Let $w$ be any vertex of $F$ such that $w \neq v$, and let $I'=(I-\{v\})\cup \{w\}$.
Note that $I'$ is independent, $|I'|=|I|$,  and yet
$|I'\cap J|=|I \cap J|-1$ contradicting our choice of $I$ and $J$.
Therefore, $G$ is well-covered.    \qed


\begin{thebibliography}{999999}




\bibitem{f2009}
A.~O.~Fradkin, On the well-coveredness of Cartesian products of graphs. {\it Discrete Math.}
{\bf 309}, 238--246 (2009).

\bibitem{h1999}
B.L.~Hartnell,  Well-covered graphs. {\it  J. Combin. Math. Combin. Comput.}
{\bf 29}, 107--115 (1999).

\bibitem{mdp1970}
M.D.~Plummer,  Some covering concepts in graphs. {\it J. Combinatorial Theory}
{\bf 8}, 91--98(1970).


\bibitem{mdp1993}
M.D.~Plummer,  Well-covered graphs: a survey.  {\it Quaestiones Math.}  {\bf 16}
 no. 3, 253--287 (1993).

\bibitem{tv1992}
J.~Topp and L.~Volkmann, On the well coveredness of products of graphs. {\it Ars Combin.}
{\bf 33}, 199--215 (1992).



\end{thebibliography}
\end{document}